\documentclass{amsart}
\usepackage{amssymb}
\newtheorem{thm}{Theorem}[section]
\newtheorem{cor}[thm]{Corollary}

\newtheorem{prop}[thm]{Proposition}
\theoremstyle{definition}

\theoremstyle{remark}

\numberwithin{equation}{section}

\newcommand{\R}{\mathbb R}
\newcommand{\Z}{\mathbb Z}
\newcommand{\eps}{\varepsilon}

\newcommand{\rt}{\rightarrow}

\newcommand{\F}{\mathcal {F}}

\newcommand{\C}{{\mathbb C}}

\begin{document}

\title[Nonnegative isotropic curvature]{Manifolds with nonnegative isotropic curvature}
\author{harish seshadri}
\address{department of mathematics,
Indian Institute of Science, Bangalore 560012, India}
\email{harish@math.iisc.ernet.in}


\thanks{This work was supported by DST Grant No. SR/S4/MS-283/05}
\begin{abstract}

We prove that if $(M^n,g)$, $n \ge 4$, is a compact, orientable,
locally irreducible Riemannian manifold with nonnegative isotropic
curvature, then one of the following possibilities hold:
\vspace{2mm}

(i) $M$ admits a metric with positive isotropic curvature

(ii) $(M,g)$ is isometric to a locally symmetric space

(iii) $(M,g)$ is K\"ahler and biholomorphic to $\C P^\frac
{n}{2}$.

(iv) $(M,g)$ is quaternionic-K\"ahler. \vspace{2mm}

This is implied by the following result: \vspace{2mm}

Let $(M^{2n},g)$ be a compact, locally irreducible K\"ahler
manifold with nonnegative isotropic curvature. Then either $M$ is
biholomorphic to $\C P^n$ or isometric to a compact Hermitian
symmetric space. This answers a question of Micallef and Wang in
the affirmative. \vspace{2mm}

The proof is based on the recent work of S. Brendle and R. Schoen
on the Ricci flow.

\end{abstract}
\maketitle
\section{Introduction}

A Riemannian manifold $(M,g)$ is said to have nonnegative
isotropic curvature (NIC) if
$$\ R_{1313} +R_{1414}+R_{2323}+R_{2424}-2R_{1234} \ge 0$$
for every orthonormal 4-frame $\{e_1,e_2,e_3,e_4\}$.

In the case of strict inequality above we say that the manifold
has positive isotropic curvature (PIC). This notion was introduced
by M. Micallef and J. Moore in ~\cite{mm} where they proved that
every compact simply-connected manifold with PIC is homeomorphic
to a sphere. In this paper, we show that the study of compact
manifolds (throughout this paper a ``compact manifold" will mean a
compact manifold without boundary) with NIC reduces to the study
of those with PIC. Our main result is the following theorem:

\begin{thm}\label{ma}
Let $(M^n,g)$, $n\ge4$, be a compact, orientable, locally
irreducible Riemannian manifold with nonnegative isotropic
curvature, then one of the following holds:

(i) $M$ admits a metric with positive isotropic curvature

(ii) $(M,g)$ is locally symmetric

(iii) $(M,g)$ is K\"ahler and biholomorphic to $\C P^{\frac
{n}{2}}$.

(iv) $(M,g)$ is quaternionic-K\"ahler.
\end{thm}

We note that according to a recent work of S. Brendle ~\cite{bre},
{\it compact Einstein manifolds with NIC are locally
symmetric.} Since quaternionic-K\"ahler manifolds are Einstein,
this result of Brendle implies that case (iv) in Theorem \ref{ma}
is included in case (iii).

A version of Theorem \ref{ma} was proved in dimension $4$ by
Micallef and Wang~\cite{mw}. \vspace{2mm}

Theorem \ref{ma} is based on the classification of K\"ahler
manifolds with NIC. We prove the following result which was
conjectured by Micallef and Wang in ~\cite{mw}.

\begin{thm}\label{mai}
Let $(M^{2n},g)$, $2n= dim_\R M \ge 4$, be a compact locally
irreducible K\"ahler manifold with nonnegative isotropic
curvature. Then either $M$ is biholomorphic to $\C P^  {n}$ or
isometric to a compact Hermitian symmetric space.
\end{thm}

Theorem \ref{mai} along with the uniqueness (up to scaling) of the
K\"ahler-Einstein metric on $\C P^n$ implies the following:

\begin{cor}\label{ke}
Let $(M,g)$, $dim_\R M \ge 4$, be a compact locally irreducible
K\"ahler-Einstein manifold with nonnegative isotropic curvature.
Then $(M,g)$ is isometric to a compact Hermitian symmetric space.
\end{cor}

{\bf Remarks:} \vspace{2mm}

(i) It follows immediately from Theorem \ref{ma} that if $(M,g)$
is a compact orientable locally irreducible Riemannian manifold
with NIC and the dimension of $M$ is odd (which is precisely the
case not treated in ~\cite{mw}), then either $M$ admits a
PIC metric or $(M,g)$ is locally symmetric. \vspace{2mm}

(ii) The four cases in Theorem \ref{ma} are not mutually
exclusive. An example of a metric satisfying cases (i) and (ii) in
Theorem \ref{ma} is one of constant positive sectional curvature.
The following remark, the proof of which is in Section 4, states
that this is the only metric for which that happens: \vspace{2mm}

{\it If $(M,g)$ is one of the following: A compact locally
symmetric space of nonconstant sectional curvature, or a  compact
K\"ahler manifold or a positive quaternionic-K\"ahler manifold,
then $M$ does not admit a metric of positive isotropic curvature}.
\vspace{2mm}

(iii) In ~\cite{mw}, Micallef and Wang  describe the structure of
{\it reducible} manifolds with NIC. Combining their result with
Theorem \ref{ma}, one can reduce the study of manifolds with NIC
to those with PIC. \vspace{2mm}


\vspace{2mm}


(iv) The following result is an easy corollary of the
Brendle-Schoen theorems mentioned below. The proof of this remark
is at the end of this paper. \vspace{2mm}

{\it Let $(M,g)$ be a compact Riemannian manifold with NIC
everywhere and PIC at some point. Then $M$ admits a metric with
PIC everywhere.}

\vspace{2mm}

(v) Results analogous to Theorem \ref{ma} are known (cf.
~\cite{bs2}, ~\cite{bw}, ~\cite{nw}) for manifolds with
nonnegative curvature operator and manifolds with non-strictly
quarter-pinched sectional curvature. Even though NIC is implied by
either of these curvature conditions, Theorem \ref{ma} does not
directly imply the earlier results. \vspace{3mm}

The proofs of Theorems \ref{ma} and \ref{mai}  are based on the
fundamental papers ~\cite{bs1}, ~\cite{bs2} of S. Brendle and R.
Schoen. According to their work, if  $g(t)$ is the solution to
Ricci flow begining at a metric $g$ with NIC, then $g(t)$ has NIC
for all $t$. Moreover, if $F_t$ denotes the set of orthonormal
4-frames on which the isotropic curvature vanishes, then for $t
>0$, $F_t$ is invariant under parallel translation by the
Levi-Civita connection of $g(t)$. From this one quickly sees,
using the Berger Holonomy Theorem, that $g(t)$ is either PIC or
$g(t)$ has holonomy in $U(m)$ or $Sp(m)Sp(1)$, i.e. $(M,g(t))$ is
K\"ahler or quaternionic-K\"ahler.  \vspace{2mm}

Suppose now that $(M,g(t))$ is K\"ahler. If $J$ is the almost
complex structure, then one knows that $F_t$ contains all
orthonormal frames of the form $\{e, J(e), f, J(f) \}$, where
$e,f$ are real tangent vectors. If these are the only elements of
$F_t$, then a version of Frankel's Conjecture due to W. Seaman
implies that $M$ is biholomorphic to $\C P^n$. On the other hand,
if $F_t$ has other elements, we will show that the holonomy is a
proper subgroup of $U(n)$. Then Berger's theorem will imply that
$(M,g(t))$ is symmetric. \vspace{2mm}

Finally, by taking a sequence $t_i \rt 0$, we can draw the same
conclusions (biholomorphic to $\C P^n$ or locally symmetric or
quaternionic-K\"ahler) about $g =lim_{i \rt \infty}g(t_i)$.

\vspace{2mm}

We remark that Ricci flow was used to study NIC and PIC in
~\cite{ham} and ~\cite{mw}. \vspace{2mm}

{\bf Acknowledgement:} The author is very grateful to Simon
Brendle for pointing out errors in the original draft  and to
Claude LeBrun for helpful suggestions.

\section{Nonnegative isotropic curvature and the Ricci flow}
In this section we recall the results of Brendle and Schoen. Let
$(M,h)$ be a compact manifold with NIC. Let $g(t), \ t \in [0,
\eps)$ be the solution to the Ricci flow equation $$\frac {
\partial g}{\partial t}=-2Ric, \ \ \ g(0)=h.$$

For  $t \in (0, \eps) $, let
$$\F_t = \bigcup_{p \in M}  \{(e_1,e_2,e_3,e_4) \ \vert \
e_i \in T_pM, \ \ g(t)(e_i,e_j)=\delta_{ij}, \ i,j=1,...,4 \}$$ be
the bundle of $g(t)$-orthonormal $4$-frames in $M$. Consider the
subset $F_t \subset \F_t$ defined by
$$F_t : = \{ (e_1,e_2,e_3,e_4) \ \vert  \
R_{1313} +R_{1414}+R_{2323}+R_{2424}-2R_{1234}=0 \}, $$ where $R$
denotes the curvature tensor of $g(t)$. i.e. $F_t$ consists of all
$g(t)$-orthonormal 4-frames $(e_1,e_2,e_3,e_4)$ at all points of
$M$ where the isotropic curvature is zero.  The two basic results
of Brendle and Schoen that we need are the following: For $t \in
(0, \eps)$, \vspace{2mm}

(i) $g(t)$ has NIC  (Section 2 in ~\cite{bs1}).

(ii) $F_t$ is invariant under parallel transport by the
Levi-Civita connection of $g(t)$ (Proposition 5 in ~\cite{bs2}).
\vspace{2mm}

\section{K\"ahler manifolds with nonnegative isotropic curvature}
This section is devoted to the proof of Theorem
\ref{mai}.\vspace{2mm}

Let $(M,g,J)$, $dim_\R M=2n$, be a compact K\"ahler manifold with
NIC with $J:TM \rt TM$ denoting the almost-complex structure. Let
$g(t)$ be the solution to Ricci flow with $g(0)=g$. We fix a $t>0$
and denote $g(t)$ by $g$ and $F_t$ by $F$. Note that by choosing
$t$ sufficiently small we can assume that $g(t)$ is locally
irreducible. Otherwise, $g=lim_{t \rt 0}g(t)$ would be locally
reducible (this can seen by considering holonomy groups).

Moreover, since $(M,h)$ is K\"ahler, so is $(M,g)$. This follows
from
\begin{thm}[R. Hamilton ~\cite{for}]
Let $(M,h)$ be a compact Riemannian $n$-manifold. If the
restricted holonomy group of $g$ lies in a subgroup of $SO(n)$
then it continues to do so for $t>0$ under the Ricci flow starting
at $h$.
\end{thm}

Note that since the restricted holonomy group $Hol^0$ of $(M,g)$
is contained in $U(n)$ by the above theorem, either $(M,g)$ is a
Hermitian symmetric space or $Hol^0=U(n)$. In the latter case,
$Hol$ is contained in the normalizer of $Hol^0$ in $SO(2n)$ (by
the orientability of $M$). This normalizer is precisely $U(n)$ and
hence $Hol=U(n)$. In either case, $(M,g)$ is K\"ahler.
\vspace{2mm}

 It can be easily checked that for any K\"ahler manifold $(M,g,J)$, 4-frames
of the form $(v,J(v),w,J(w)) \in F$ , for any $v,w \in T_pM$
satisfying $g(v,w)= g(v,Jw)= 0$.

The analysis of the Ricci flow splits into two cases, depending on
the set $F$.

\vspace{2mm}

{\bf Case I}: $F$ does not contain any frame of the form
$(u,J(u),J(v),v)$. \vspace{2mm}

We claim that $(M,g)$ is biholomorphic to $\C P^n$ in this case.
\vspace{2mm}

For any $p \in M$, let $u,v \in T_pM$ with $g(u,v)=g(u,J(v))=0$.
As mentioned above, by the symmetries of the curvature tensor of a
K\"ahler manifold, we have $(u,J(u),v,J(v)) \in F$, i.e.,
\begin{align}\label{one}
&0 \ = \  R(u,v,u,v)+R(u,J(v),u,J(v))+R(J(u),v,J(u),v) \\ \notag
& \ \ \ \ \ \ +R(J(u),J(v),J(u),J(v))-2R(u,J(u),v,J(v)).\notag \\
&= \ 2 \Bigl ( R(u,v,u,v)+R (u,J(v),u,J(v) - R(u,J(u),v,J(v)).
\Bigr ) \notag
\end{align}

Since $(u,J(u), J(v),v) \notin F$, we have
\begin{align}\label{two}
&0 \ < \  R(u,J(v),u,J(v))+R(u,v,u,v)+R(J(u),J(v),J(u),J(v)) \\
\notag
& \ \ \ \ \ \ \ +R(J(u),v,J(u),v)-2R(u,J(u),J(v),v). \notag \\
&= \ 2 \Bigl ( R(u,v,u,v)+R (u,J(v),u,J(v) + R(u,J(u),v,J(v)).
\Bigr ) \notag
\end{align}
Adding (~\ref{one}) and (~\ref{two}), we get
\begin{equation}\label{sat}
 R(u,v,u,v) \ + \ R
(u,J(v),u,J(v)) \ > 0
\end{equation}
for every $p \in M$ and $u, v \in T_pM$ with $g(u,v)=g(u,J(v))=0$.
This condition, sometimes referred to as {\it orthogonal
bisectional curvature} in the literature, is precisely the
condition $(\ast \ast)$ on Page 846 of ~\cite{sea}. \vspace{2mm}

We now state a version of the Frankel Conjecture, following
~\cite{sea} and ~\cite{mw}: \vspace{2mm}

{\it If $(M^{2n},g)$ is a compact locally irreducible K\"ahler
manifold with nonnegative isotropic curvature satisfying
(\ref{sat}), then $M$ is biholomorphic to $\C P^n$.} \vspace{2mm}

{\it Proof}: The proof is the same as in the paper of W. Seaman,
~\cite{sea}, except for some minor changes. For the sake of
completeness we outline it, emphasizing the parts where the
curvature assumptions are used. Note that the sign convention for
the curvature tensor in ~\cite{sea} and ~\cite{mw} is the opposite
of what we follow. Seaman's version of the Frankel Conjecture,
Theorem B of ~\cite{sea}, asserts the following: \vspace{2mm}

{\it If $(M^{2n},g)$ is a compact K\"ahler manifold with
nonnegative isotropic curvature satisfying (\ref{nag} ), then $M$
is biholomorphic to $\C P^n$} \vspace{2mm}

where (\ref{nag}) is the condition that
\begin{equation} \label{nag}
 R_{ikik}+ R_{ilil}+R_{jkjk}+R_{jljl} \ > \ 0
\end{equation}

for all orthonormal vectors $e_i,e_j,e_k,e_l$. Our main
observation, which is justified in the steps below, is that we can
replace (\ref{nag}) by (\ref{sat}) and local irreducibility of the
metric to obtain the same conclusion.

Seaman's proof closely follows the proof of the Frankel Conjecture
for positive bisectional curvature by Siu and Yau ~\cite{sy}.
There are three points where the curvature assumptions play a
role: \vspace{2mm}

{\bf (i)} To start the proof one needs to know that $\pi_2(M)=
{\mathbb Z} \oplus \ torsion $. This is done by showing that
$b_2(M)=1$ and $\pi_1(M)=\{ 0 \}$ and applying the Hurewicz
theorem. A reading of the proof by Seaman shows that (\ref{nag})
is used to prove that $b_2(M)=1$. In our case, we do not have
(\ref{nag}) but instead have local irreducibility of the metric.
This will suffice since Theorem 2.1 (b) of Micallef-Wang
~\cite{mw} (note that $b_2(M) \neq 0$ since $M$ is K\"ahler)
implies: {\it Let $(M^{2n},g)$ be a compact, locally irreducible
Riemannian manifold with NIC and $b_2(M) \neq 0$. Then $b_2(M)=1$
and $M$ is simply-connected.}

As will be seen below {\it the only curvature assumptions needed
in steps (ii) and (iii) are nonnegative isotropic curvature and
positive orthogonal bisectional curvature (\ref{sat})}.

\vspace{2mm}

{\bf (ii)} Choose a generator $a$ of $H^2(M, \Z)= \Z$ which is a
negative multiple of the K\"ahler class and $\alpha$ be a
generator of the free part of $\pi_2(M)$ such that $a(\alpha)=1$.
By a result of Kobayashi-Ochiai ~\cite{ko} to show that $M$ is
biholomorphic to $\C P^n$, it is enough to show that
$c_1(M)(\alpha) \ge n+1$, where $c_1(M)$ denotes the first Chern
class of $M$. By the Sacks-Uhlenbeck theorem (~\cite{su}, Theorem
5.5) we may represent the free homotopy class of $\alpha$ by
$\sum_{i=1}^k f_i$ where each $f_i$ is harmonic and energy
minimizing in its free homotopy class. One claims that each $f_i$
is $\pm$ holomorphic. To see this, one needs the complex
formulation of isotropic curvature and the second variation
formula for the energy of minimal surfaces as in the paper of
Micallef-Moore ~\cite{mm}. We recall this below and prove the
holomorphicity lemma in detail, following Seaman:

Let $(M,g)$ be a Riemannian manifold. Let $E=TM \otimes \C$.
Extend $g$ to a symmetric bilinear form $(,)$ and a Hermitian form
$\langle \langle , \rangle \rangle$ on $E$. For $p \in M$, an
element $v \in E_p$ is said to be {\it isotropic} if $(v,v)=0$. A
$2$-plane $P \subset E_p$ is isotropic if every element of $P$ is
isotropic. It can be checked that $\{ v,w \}$ spans an isotropic
$2$-plane if and only if $v$ and $w$ are linearly independent and
$$(v,v) \ = \ (w,w) \ = \ (v,w) \ = 0.$$
Also, $Span_\C \{ v, w \}$ is an isotropic $2$-plane if and only
if there exist orthonormal vectors $e_1,e_2,e_3,e_4$ in $T_pM$
such that
\begin{equation} \notag
 v \ = \ \frac { \Vert v \Vert }{\sqrt 2}
\Bigl ( e_1 + \sqrt {-1} e_2 \Bigr ), \ \ \ \ w \ = \ \frac {
\Vert w \Vert }{\sqrt 2} \Bigl ( e_3 + \sqrt {-1} e_4 \Bigr ),
\end{equation}
where the norm is with respect to $\langle \langle , \rangle
\rangle$.

 \vspace{2mm}

Let
$${\mathcal R}: {\bigwedge} ^2 E_p \rt {\bigwedge} ^2 E_p $$
denote the complex linear extension of the curvature operator. If
$v,w, e_1,..,e_4$ are as above, then it can be checked that
\begin{equation} \label{koth}
4 \langle \langle \ {\mathcal R} ( v \wedge w), \ v \wedge w \
\rangle \rangle \ =  \Vert v \Vert^2 \Vert w \Vert^2 \
R^{iso}(e_1,e_2,e_3,e_4)
\end{equation} \vspace{2mm}
where for any tangent vectors $X,Y,Z,W$
\begin{align}\notag
R^{iso}(X,Y,Z,W) :=&R(X,Z,X,Z)+R(X,W,X,W)+R(Y,Z,Y,Z) \\ \notag
&+R(Y,W,Y,W)-2R(X,Y,Z,W). \notag
\end{align}
The assumption in Case I is that
$$R^{iso}(X,J(X),J(Y),Y) >0$$ for any orthonormal frame $(X,J(X),J(Y),Y)$.

Let $f: S^2 \rt M$ be a smooth immersion. Consider $F = f^\ast (E)
= f^\ast TM \otimes \C$. We pull-back $g, (  ,  )$ and $\langle
\langle, \rangle \rangle$ and denote them by the same symbols. $F$
also carries the pull-back of the Levi-Civita connection on $M$,
extended complex linearly. This connection is Hermitian with
respect to $\langle \langle, \rangle \rangle$. It is well-known
that a Hermitian bundle $V$ on a Riemann surface $\Sigma$ with a
Hermitian connection $\nabla$ can be endowed with a holomorphic
structure $\overline
\partial$ in which a section $s$ of $V$ is holomorphic is and
only if
$$ \nabla_{\frac {\partial }{\partial \bar z}}s = 0$$
in any local holomorphic coordinate $z$ on $\Sigma$. \vspace{3mm}

Consider $S^2$ as $\C \cup \{ \infty \}$ and let $Z$ be the
holomorphic vector field which is $\frac {\partial} {\partial z}$
on $\C$ and $0$ at $\infty$. If $f$ is a harmonic map, then
$f_\ast Z$ is a holomorphic isotropic section of $F$, where
$f_\ast: TS^2 \otimes \C \rt F$. If, in addition, $f$ is stable
and $(M,g)$ has nonnegative isotropic curvature then it follows
from the stability inequality of Micallef and Moore that
\begin{equation}\label{kk}
\langle \langle \ {\mathcal R} ( s \wedge f_\ast Z), \ s \wedge
f_\ast Z \ \rangle \rangle \ = \ 0,
\end{equation}
for any holomorphic section $s$ of $F$ such that $s$ and $f_\ast
Z$ span an isotropic $2$-plane. \vspace{2mm}

Now we can prove that {\it a harmonic map into a K\"ahler manifold
with nonnegative isotropic curvature and positive orthogonal
bisectional curvature has to be $\pm$ holomorphic i.e., $J \circ
f_\ast (Z) = \pm \sqrt {-1} f_\ast (Z)$.}

\begin{proof}
Let $S= \{ x \in S^2 \ \vert \ (f_\ast Z)_x =0 \}$. Then $S$ is a
finite set. Take $p \in S^2 \setminus S$.


Let $w =(f_\ast Z)_p$ and $F_p = f^\ast T_pM \otimes \C$ as
before. We can write $F_p=F_p^{(1,0)} \oplus F_p^{(0,1)}$, where
$F_p^{(1,0)}$ , \ $F_p^{(0,1)}$ are the eigenspaces of $J$ for the
eigenvalues $\sqrt {-1}, \ -\sqrt {-1}$ respectively.

It is enough to prove that $w$ is an eigenvector of $J$. Suppose
not. Now $w$ and $J(w)$ span an isotropic $2$-plane. Also, we can
write $w=w' + w''$ as the sum of non-zero $(1,0)$ and $(0,1)$
parts. Since $w$ is also isotropic, we can find orthonormal $e_1,
e_2 \in f^\ast T_pM$ such that $g( e_1, J(e_2)) =0$ and
$$ \sqrt 2 \frac {w'} {\Vert w' \Vert} = e_1 - \sqrt {-1} \ J(e_1), \
\ \sqrt 2 \frac {w''} {\Vert w'' \Vert} = e_2 + \sqrt {-1} \
J(e_2).$$ Since $w \wedge J(w)  = -2 \sqrt {-1} \ w' \wedge w''$,
(\ref{koth}) gives
\begin{align} \notag
& \langle \langle \ {\mathcal R} ( w \wedge J(w)) , \ w \wedge
J(w) \ \rangle \rangle \ = \ 4 \langle \langle \ {\mathcal R} ( w'
\wedge w'') , \ w' \wedge w'' \ \rangle \rangle \notag \\
&= \Vert w' \Vert^2  \Vert w'' \Vert^2 (R^{iso}(e_1, -J(e_1),e_2,
J(e_2))) \notag \\
&= \Vert w' \Vert^2  \Vert w'' \Vert^2 (R^{iso}(e_1,
J(e_1),J(e_2), e_2)) >0, \notag
\end{align}

\noindent which contradicts (\ref{kk}).
Hence $w$ is
an eigenvector for $J$ on the complement of finitely many points
in $S^2$ i.e. $f$ is $\pm$ holomorphic on this set. It follows
that $f$ is $\pm$ holomorphic on all of $S^2$. \vspace{2mm}

{\bf (iii)} By (ii) we have

\begin{equation}\label{qx}
\alpha = \sum_{i=1}^k f_i,
\end{equation}

\noindent where each $f_i$ is $\pm$ holomorphic. On pages 853-854
of ~\cite{sea} Seaman proves the following fact:

{\it Let $(M^n,g)$ be a compact K\"ahler manifold with positive
orthogonal bisectional curvature and let $TM$ denote the
holomorphic tangent bundle of $M$. Let $f: S^2 \rt M$ be any
nonconstant holomorphic map. By a theorem of Grothendieck ,
$f^\ast TM =L_1 \oplus ...\oplus L_n$ where the $L_i$ are
holomorphic line bundles. The claim is that}
$$ c_1(L_i) \ge 1  \ \ \forall \ i \ \ \ {\rm and} \ \ c_1(L_{i_0}) \ge
2 \ \ {\rm for \ some} \ i_0.$$

We refer the reader to ~\cite{sea} for a proof.

This result implies that if $k=1$ in (\ref{qx}) and $f_1$ is
holomorphic then $c_1(M)(\alpha) =c_1(M)[f_1^\ast(S^2)]
=c_1(f_1^\ast(TM))[S^2] \ge n+1$. Actually $k=1$ implies that
$f_1$ is holomorphic since $c_1(M)$ is a negative multiple of the
K\"ahler class. In order to show that $k=1$ one proves the
Deformation Lemma (Proposition 3, page 201 of ~\cite{sy}). Given
this Lemma the proof that $k=1$ is identical to that of Claim 2 of
~\cite{sy}, pages 202-203.

As was observed by Futaki ~\cite{fut}, the Deformation Lemma holds
for any compact K\"ahler manifold and any holomorphic map $f:S^2
\rt M$ as long as $f^\ast TM$ splits as a direct sum of positive
line bundles. See ~\cite{fut} or the last part of Page 854 of
~\cite{sea} for a proof.

\end{proof}

\vspace{2mm}

Hence we conclude that $M$ is biholomorphic to $\C P^n$ under the
assumptions of Case I.

\vspace{2mm}

{\bf Case II:} $F$ contains a frame of the form $(u,J(u),
J(v),v)$.

\vspace{2mm}

We claim that the restricted holonomy group of $(M,g)$ cannot be
the whole group $U(n) \subset SO(2n)$. This will complete the
proof of Theorem \ref{mai} by Berger's Holonomy Theorem.
\vspace{2mm}

Suppose $Hol^0(M)=U(n)$. \vspace{2mm}

The fact that $(u,J(u),v,J(v))$ and $(u,J(u),J(v),v)$ both belong
to $F$ gives (as in (\ref{one}) and (\ref{two}))
\begin{align}\notag
R(u,J(u),v,J(v)) &= R(u,J(v),u,J(v))+R(u,v,u,v) \\ \notag &
+R(J(u),J(v),J(u),J(v))+R(J(u),v,J(u),v)=0 \notag
\end{align}
which, by the symmetries of K\"ahler curvature, gives
\begin{equation} \label{ah}
R(u,J(u),v,J(v))=R(u,J(v),u,J(v))+R(u,v,u,v)=0
\end{equation}

We now use the holonomy action of $U(n)$ repeatedly to get a
contradiction. First, there is an element of $U(n)$ under which
\begin{equation} \label{lat1}
u \rt \frac {1}{\sqrt 2}(u -J(v)), \ \ v \rt \frac {1}{\sqrt 2}(u
+ J(v)).
\end{equation}

Hence $(\frac {1}{\sqrt 2}(u -J(v)),  \frac {1}{\sqrt 2}(J(u) +v),
\frac {1}{\sqrt 2}(u +J(v)), \frac {1}{\sqrt 2}(J(u) -v)) \in F$.

 The equation corresponding to (\ref{ah}) is
$$R(\frac {1}{\sqrt 2}(u -J(v)), \ \frac {1}{\sqrt 2}(J(u) +v), \frac {1}{\sqrt 2}(u +J(v)), \frac
{1}{\sqrt 2}(J(u) -v))=0.$$

This equation gives
\begin{align} \notag
&R(u-J(v), \ J(u)+v,  \ u+J(v), \ J(u)-v) \notag \\
& = \  R(u,\ J(u), \ u, \ J(u))  +  R(u,\ J(u), \ u, \ -v)  +  R(u, \ J(u), \ J(v), \ J(u)) \notag \\
& + \ R(u, \ J(u), \ J(v), \ -v)  +  R(u, \ v, \ u, \ J(u))  +  R(u, \ v, \ u, \ -v) \notag \\
& + \  R(u, \ v, \ J(v), \ J(u))  +  R(u, \ v, \ J(v), \ -v)  +  R(-J(v), \ J(u), \ u, \ J(u)) \notag \\
& + \ R(-J(v), \ J(u), \ u, \ -v)  +   R(-J(v), \ J(u), \ J(v), \ J(u))  +  R(-J(v), \ J(u), \ J(v), \ -v) \notag \\
& + \ R(-J(v), \ v, \ u, \ J(u))  +  R(-J(v), \ v, \ u , \ -v)  +  R(-J(v), \ v, \ J(v), \ J(u)) \notag \\
& + \ R(-J(v), \ v, \ J(v), \ -v) \ = \ 0. \notag
\end{align} \vspace{2mm}

The sum of the 2nd, 3rd, 5th and 9th terms on the right side is
$$R(u,J(u),u,-v) \ + \ R(u,  J(u),  J(v),  J(u)) \ + \ R(u,  v,  u,
J(u)) \ + \ R(-J(v),  J(u), u, J(u)) =0.$$

Similarly the sum of the 8th, 12th, 14th and 15th terms is
$$R(u,v,J(v),-v) \ + \ R(-J(v),J(u),J(v),-v) \ + \  R(-J(v),v,u ,-v) \  +  \
R(-J(v),v,J(v),J(u))=0.$$

The sum of the 6th, 7th, 10th and 11th terms is
\begin{align}\notag
& R(u, \ v, \ u, \ -v)+ \  R(u, \ v, \ J(v), \ J(u)) \\ \notag &
+R(-J(v), J(u), \ u, \ -v) +   R(-J(v), \ J(u), \ J(v), \ J(u)) =
 -4 R(u, \ v, \ u, \ v). \notag
\end{align}

The sum of the 4th and 13th terms is
$$  R(u, \ J(u), \ J(v), \ -v) \ + \ \ R(-J(v), \ v, \ u, \ J(u))
=0$$ by (\ref{ah}).

Using these four equations and simplifying the expansion of
$R(u-J(v), \ J(u)+v,  \ u+J(v), \ J(u)-v)$, we get \vspace{2mm}
\begin{align}
R(u,J(u),u,J(u)) \ + \ R(v,J(v),v,J(v)) \ - \ 4R(u,v,u,v)=0.
\label{eh}
\end{align} \vspace{2mm}

Next, consider the element of $U(n)$ which takes
\begin{equation}
 u \rt \frac {1}{\sqrt 2}(u -v), \ \ v \rt \frac {1}{\sqrt 2}(u +v).
\end{equation}
Now the corresponding equation corresponding to (\ref{eh})
(obtained by substituting $-J(v)$ for $v$ in (\ref{eh})) is
\vspace{2mm}

\begin{align}
R(u,J(u),u,J(u)) \ + \ R(v,J(v),v,J(v)) \ - \ 4R(u,J(v),u,J(v))=0.
\label{eh2}
\end{align} \vspace{2mm}

Combining (\ref{eh2}), \ (\ref{eh}) and (\ref{ah}) gives
\begin{align}\notag
&R(u,v,u,v)=R(u,J(v),u,J(v)) \notag \\
&= R(u,J(u),u,J(u))+R(v,J(v),v,J(v))=0. \label{eh3}
\end{align}

Extend $\{ e_1=u, e_2=v \}$ to an orthonormal basis
$\{e_1,J(e_1),e_2,J(e_2),...,e_n,J(e_n) \}$. By considering the
element of $U(n)$ which interchanges $u$ and $e_i$, $J(u)$ and
$J(e_i)$ ($i \ge 3$) and keeps the other elements of the basis
fixed, we see that $(e_i,J(e_i),J(v),v) \in F$. A similar
operation on $v$ shows that $(e_i, J(e_i),e_j,J(e_j)) \in F$ for
all $1 \le i \le n$, $1 \le j \le n$, $i \neq j$. The equations
corresponding to (\ref{eh3}) are
\begin{align}\label{eh4}
&R(e_i,e_j,e_i,e_j)=R(e_i,J(e_j),e_i,J(e_j)) \\ \notag
&=R(e_i,J(e_i),e_i,J(e_i))+R(e_j,J(e_j),e_j,J(e_j))=0. \notag
\end{align}

The equations
$R(e_i,J(e_i),e_i,J(e_i))+R(e_j,J(e_j),e_j,J(e_j))=0$ for all $i
\neq j$ clearly imply that $R(e_i,J(e_i),e_i,J(e_i))=0$ for all
$i$. In particular $R(u,J(u),u,J(u))=0$ . If $w$ is an arbitrary
unit vector, then we can find an element $T$ in $U(n)$ taking $u$
to $w$. By considering the frame $(w, J(w), J(T(v)), T(v)) \in F$
we get $R(w,J(w),w,J(w))=0$. Hence the holomorphic sectional
curvature of $g$ is zero which implies that $g$ is flat and hence
locally reducible.




This contradiction completes the proof of Theorem \ref{mai}.
\hfill $\square$ \vspace{2mm}

\section{The classification theorem}
We now outline the proof of Theorem \ref{two}: Let $(M,g)$ be a
compact, orientable, locally irreducible  manifold of nonnegative
isotropic curvature. Let $g(t), \ t \in [0, \delta) $ denote the
solution to Ricci flow with $g(0)=g$. By ~\cite{bs1}, $g(t)$ has
nonnegative isotropic curvature. Suppose $g(t)$ does not have
strictly positive isotropic curvature for any $t \in (0,\delta)$.
Then $F_t \neq \phi$ for all $t \in (0,\delta)$. This implies that
$Hol^0(M,g(t)) \neq SO(n)$ for all such $t$: Suppose
$Hol^0(M,g(t_0))=SO(n)$ for some $t_0$. This assumption along with
the invariance of $F_{t_0}$ under parallel transport implies that
every orthonormal 4-frame at every point would be in $F_{t_0}$,
i.e., every isotropic curvature at every point would be zero. As
scalar curvature can be expressed as the sum of isotropic
curvatures (see Proposition \ref{dd} below), this implies that
$R=0$ for $(M,g(t_0))$. Since the initial metric $g$ has
nonnegative scalar curvature, by the maximum principle for scalar
curvature along the flow, the scalar curvature of $(M,g(t))$ is
identically zero for all $t \in [0, t_0]$. From the evolution
equation for $R$ under Ricci flow, $$ \frac {\partial R}{\partial
t} = \triangle R + \vert Ric \vert^2,$$ it follows that
$Ric(g)=0$.

Moreover we have the following \vspace{2mm}

\begin{prop}[Proposition 2.5 of ~\cite{mw}] ~\label{dd}
A metric $g$ with positive (resp. nonnegative) isotropic curvature
has positive (resp. nonnegative) scalar curvature. If the scalar
curvature is identically zero then $g$ must be conformally flat.
\end{prop}

Therefore $g$ is also conformally-flat and hence flat,
contradicting local irreducibility.

Hence we conclude that $Hol^0 (M,g(t))$ is a proper subgroup of
$SO(n)$ for all $t \in (0,\delta)$. As before, for sufficiently
small $t$, say $t < \delta_1$, $g(t)$ will also be locally
irreducible.

Take any $t' \in (0,\delta_1)$. If $(M,g(t'))$ is not locally
symmetric, then $Hol^0$ would have to be either $U(m)$ or
$Sp(m)Sp(1)$, with $m$ equal to $\frac {n}{2}$ or $\frac {n}{4}$,
respectively. The other possibilities for $Hol^0$ in the Berger
Holonomy theorem can be ruled out since they would again lead to
scalar flatness or local symmetry.

In case $Hol^0=U(m)$, the full holonomy group  $Hol= U(m)$ as
well, since $Hol$ is contained in the normalizer of $Hol^0$ in
$SO(n)$ (note that we are assuming the orientability of $M$ here).
This normalizer is precisely $U(n)$.
Since $Hol=U(m)$, $(M,g(t'))$ is K\"ahler and Theorem \ref{mai}
implies either $(M,g)$ is biholomorphic to $\C P^n$ or a symmetric
space.

If $Hol^0=Sp(m)Sp(1)$, then again $Hol= Sp(m)Sp(1)$ (note that
since $Sp(m)Sp(1)$ is a maximal subgroup of $SO(4n)$ ~\cite{gray},
its normalizer in $SO(4n)$ is itself).
Hence for any $t' \in (0, \delta_1)$ either $(M,g(t'))$ is
K\"ahler and biholomorphic to $\C P^n$ or a locally symmetric
space or quaternionic-K\"ahler. By taking a sequence $t_i \rt 0$
as $i \rt \infty$ and noting that $g(t_i) \rt g$ in $C^\infty$, we
can conclude the same about $(M,g)$. This completes the proof of
Theorem \ref{ma}.

 \hfill $\square$ \vspace{2mm}

{\it Proof of Remark (ii) of Page 2}: First, by the results of
~\cite{mw} or ~\cite{sea}, if $(M,h)$ has positive isotropic
curvature, then $H^2(M,\R)=\{0\}$.  Hence $M$ cannot admit a
K\"ahler metric.  If $M$ admitted a positive quaternionic-K\"ahler
metric, then it would be simply-connected and $H^4(M,\R) \neq
\{0\}$. But since $M$ admits a PIC metric, it would be
homeomorphic to a sphere, a contradiction.

Suppose that there is a locally symmetric metric $g$ on $M$. Note
that $M$ admits a metric of positive scalar curvature (namely
$h$), $(M,g)$ will have to be a locally symmetric space of compact
type. In particular $(M,g)$ is Einstein with positive scalar
curvature. By the Bonnet-Myers theorem $M$ has finite fundamental
group and hence the universal cover $\tilde M$ is a compact
manifold. The Micallef-Moore theorem applied to $(\tilde M, \tilde
h)$ implies that $\tilde M$ is homeomorphic to a sphere. For
topological reasons the only locally symmetric metric on a
topological sphere is one of constant positive sectional
curvature. Therefore $\tilde g$ (and $g$) would have to be of
constant positive sectional curvature. \hfill $\square$
\vspace{2mm}

{\it Proof of Remark (iv) on Page 2}: Suppose $K^{iso} \ge 0$ on
$M$ and $K^{iso}(p) >0$ for some $p \in M$. With notation as
earlier, if $F_t$ is not empty, i.e. if $g(t)$ does not have
$K^{iso} >0$, then by the invariance under parallel transport of
$F_t$, there is a $4$-frame for which $K^{iso}=0$ at {\it every}
point of $M$. In particular there is such a frame at $p$. Hence we
have a time dependent $4$-frame $(e_1(t),e_2(t),e_3(t),e_4(t))$ at
$p$ for which $K^{iso}=0$. We can choose a sequence of times $t_i
\rt 0$ as $i \rt \infty$ for which the corresponding sequence of
frames converges to an orthonormal frame on $(M,g)$. This frame
will satisfy $K^{iso} =0$ contradicting $K^{iso} >0$ at $p$.
\hfill $\square$.

\end{document}